\documentclass[12pt, a4paper]{amsart}
\usepackage{amsmath}
\usepackage{geometry,amsthm,graphics,tabularx,amssymb,mathrsfs,shapepar}
\usepackage{amscd}
\usepackage[all]{xypic}

\newcommand{\BC}{{\mathbb {C}}}

\newcommand{\CH}{{\mathcal {H}}}

\newcommand{\CR}{{\mathcal {R}}}
\newcommand{\CS}{{\mathcal {S}}}

\newcommand{\End}{{\mathrm{End}}}

\newcommand{\GL}{{\mathrm{GL}}}

\newcommand{\Hom}{{\mathrm{Hom}}}

\newcommand{\Sp}{{\mathrm{Sp}}}

\newcommand{\cover}[1]{\widetilde{#1}}

\newcommand{\Irr}{\operatorname{Irr}}

\newcommand{\oL}{\operatorname{L}}

\newcommand{\oO}{\operatorname{O}}

\newcommand{\oU}{\operatorname{U}}

\newcommand{\g}{\mathfrak g}

\newcommand{\p}{\mathfrak p}

\renewcommand{\l}{\mathfrak l}

\newcommand{\s}{\mathfrak s}

\renewcommand{\sl}{\mathfrak s \mathfrak l}

\newcommand{\C}{\mathbb{C}}
\newcommand{\R}{\mathbb R}

\newcommand{\rF}{\mathrm{F}}

\newcommand{\rH}{\mathrm{H}}

\newcommand{\be}{\begin {equation}}
\newcommand{\ee}{\end {equation}}
\newcommand{\bee}{\begin {equation*}}
\newcommand{\eee}{\end {equation*}}

\theoremstyle{Theorem}

\theoremstyle{Theorem}

\newtheorem{thm}{Theorem}[section]

\theoremstyle{Theorem}

\theoremstyle{Theorem}
\newtheorem{prp}{Proposition}[section]

\newtheorem{conj}[thm]{Conjecture}

\theoremstyle{Plain}
\newtheorem*{remark}{Remark}

\theoremstyle{Definition}

\theoremstyle{remark}
\newtheorem{example}{Example}

\title{Local theta correspondence: the basic theory}

\author{Binyong Sun}
\address{Academy of Mathematics and Systems Science, Chinese Academy of Sciences \\
and School of Mathematical Sciences, University of Chinese Academy of Sciences, Beijing 100190, China} \email{sun@math.ac.cn}

\author{Chen-Bo Zhu}
\address{Department of Mathematics\\
National University of Singapore\\
Block S17, 10 Lower Kent Ridge Road, Singapore 119076}
\email{matzhucb@nus.edu.sg}

\begin{document}


\begin{abstract}
We give an elementary introduction to Classical Invariant Theory and its modern extension ``Transcending Classical Invariant Theory", commonly known as the theory of local theta correspondence. We explain the two fundamental assertions of the theory: the Howe duality conjecture and the Kudla-Rallis conservation relation conjecture. We give a status report on the problem of explicitly describing local theta correspondence in terms of Langlands-Vogan parameters. We conclude with a discussion on a certain problem of automatic continuity, which manifests unity of the theory in algebraic and smooth settings.
\end{abstract}

\maketitle


\section{Classical invariant theory}
Classical invariant theory is the study of the polynomial invariants for an arbitrary number of (contravariant or covariant) variables for a standard classical group action.
We will focus on the case of the compact orthogonal group. The basic references are \cite{Ho2} and \cite{KV}.

So let $G=\oO(n)$, the compact orthogonal group in $n$ variables. It has the standard action
\[
\oO(n)\curvearrowright \R^n,
\]
and it induces a representation
\[
   \oO(n)\curvearrowright \C[\R^n]
\]
given by
\[
(g\cdot f)(x):=f(g^{-1}x)\quad \textrm{for all } g\in \oO(n),\, f\in \C[\R^n], \, x\in \R^n,
\]
where $\C[\R^n]$ denotes the space of complex valued polynomial functions on $\R^n$.

Denote by $\{x_i\}_{1\leq i\leq n}$ the standard coordinates of $\R ^n$. It is elementary to show that the space of $\oO(n)$-invariant polynomials is
\[
  \C[\R^n]^{\oO(n)}=\C[q],
\]
where
\[
 q:=\sum_{i=1}^n x_i^2
 \]
 is the quadratic polynomial on $\R^n$ preserved by $\oO(n)$.

It is a part of classical invariant theory to answer the following question:
\[
  \textrm{How to decompose the representation $\oO(n)\curvearrowright \C[\R^n]$ ?}
\]
To do this, it is profitable to introduce a Lie algebra action
\be\label{actsl2}
\sl_2(\C)\curvearrowright \C[\R^n],
\ee
specified by
\[
  \begin{array}{rcl}
    \mathbf e:= \left[
                   \begin{array}{cc} 0&1\\ 0&0\\
                   \end{array}
  \right]
  &\mapsto&   -\frac{1}{2}\sum_{i=1}^n x_i^2, \smallskip\\
   \mathbf h:= \left[
                   \begin{array}{cc} 1&0\\ 0&-1\\
                   \end{array}
  \right]
  &\mapsto & \frac{n}{2}+\sum_{i=1}^n x_i\frac{\partial}{\partial x_i},\smallskip \\
    \mathbf f:=\left[
                   \begin{array}{cc} 0&0\\ 1&0\\
                   \end{array}
  \right]
  &\mapsto & \frac{1}{2} \sum_{i=1}^n \frac{\partial^2}{\partial x_i^2}.
 \end{array}
\]
The most important feature of this Lie algebra action is that it commutes with the action of $\oO(n)$, and the resulting extended action of
$\oO(n)\times \sl_2(\C)$ will help us in decomposing the action of $\oO(n)$ on $\C[\R^n]$.

The group $\oO(n)$ acts on the space of harmonic polynomials:
\be\label{on}
 \oO(n)\curvearrowright \CH[\R^n]:=\{\varphi\in \C[\R^n]\mid \mathbf f\cdot \varphi=0\}.
\ee
According to the action of $\mathbf h$, we have the obvious decomposition:
\[\CH[\R^n]=\bigoplus_{d=0}^\infty \CH^d[\R^n],
\]
where
\[
  \CH^d[\R^n]:=\{\varphi\in \mathcal H[\R^n]\mid\textrm{$\varphi$ is homogeneous of degree $d$}\}.
\]

The following result is well-known, and commonly called the theory of spherical harmonics.

\begin{thm}\label{sh1}
The representation
\[
 \oO(n)\curvearrowright \CH^d[\R^n]
 \]
  is either irreducible or zero.
Moreover, we have the following decomposition as an $\oO(n)\times \sl_2(\C)$-module:
\[
   \C[\R^n]\cong\bigoplus_{d\geq 0, \,\CH^d[\R^n]\neq 0} \CH^d[\R^n]\otimes \oL(d+\frac{n}{2}).
\]
Here $\oL(d+\frac{n}{2})$ is the irreducible lowest weight module of $\sl_2(\C)$ with the lowest weight $d+\frac{n}{2}$.
\end{thm}

More generally, we consider the action of $\oO(n)$ on the direct sum of $k$-copies of $\R^n$:
\[
\oO(n)\curvearrowright \R^{n\times k}, \qquad (k\geq 0, \textrm{ $\R^{n\times k}$ denotes the space of $n\times k$ real matrices}).
\]
Again it induces a representation
\be\label{acto}
   \oO(n)\curvearrowright \C[\R^{n\times k}].
\ee

Generalizing the action of $\oO(n)\times \sl_2(\C)$ for $k=1$, we now have a natural extension of the representation \eqref{acto}
 to a representation
\begin{equation}\label{actsp}
 \oO(n)\times \s\p_{2k}(\C)\curvearrowright \C[\R^{n\times k}].
\end{equation}
(The action of the symplectic Lie algebra $\s\p_{2k}(\C)$ may be viewed as the hidden symmetry of $\C[\R^{n\times k}]$.)
The symplectic Lie algebra $\s\p_{2k}(\C)$ has the usual decomposition
\[
  \s\p_{2k}(\C)=\p^+\oplus \g\l_k(\C)\oplus \p^-,
\]
and in the standard coordinates $\{x_{i,j}\}_{1\leq i\leq n, 1\leq j\leq k}$ of $\R ^{n\times k}$, the action is described by the following:
\[
  \begin{array}{rcl}
    \p^+  &\mapsto&   \text{span}\{\sum_{l=1}^n x_{li}x_{lj}\}_{1\leq i,j\leq k}, \smallskip\\
   \g\l_k(\C) &\mapsto & \text{span}\{\frac{n}{2}\delta_{i,j}+\sum_{l=1}^n x_{li}\frac{\partial}{\partial x_{lj}}\}_{1\leq i,j\leq k},\smallskip \\
    \p^-  &\mapsto & \text{span}\{\sum_{l=1}^n \frac{\partial^2}{\partial x_{li}\partial x_{lj}}\}_{1\leq i,j\leq k}.
 \end{array}
\]

We now have a representation
\be\label{on2}
 \oO(n)\otimes \g\l_k(\C)\curvearrowright \CH[\R^{n\times k}]:=\{\varphi\in \C[\R^{n\times k}]\mid X\cdot \varphi=0\quad \textrm{for all $X\in \p^-$}\}.
\ee

Generalizing Theorem \ref{sh1}, the following theorem is one of the main results in classical invariant theory. See \cite{Ho2} or \cite{KV}.

\begin{thm}\label{sh2}
As an $\oO(n)\times \g\l_k(\C)$-module,
\[
  \CH[\R^{n\times k}]=\bigoplus \pi\otimes \sigma
 \]
 is multiplicity free, and $\pi$ and $\sigma$ determines each other, where $\pi$ runs over a set of irreducible finite-dimensional representations of $\oO(n)$, and $\sigma$ runs over a set of irreducible finite-dimensional representations of $\g\l_k(\BC)$. As an $\oO(n)\times \s\p_{2k}(\C)$-module,
\[
  \C[\R^{n\times k}]=\bigoplus \pi\otimes \oL(\sigma),
 \]
 where $\oL(\sigma)$ is the unique irreducible quotient of the generalized Verma module $\oU(\s\p_{2k}(\C))\otimes_{\oU(\g\l_k(\C)\oplus \p^-)} \sigma$, and ``$\oU$" indicates the universal enveloping algebra.
\end{thm}

A natural question is about the occurrence, namely, for a given irreducible representation $\pi$ of $\oO(n)$, we want to know whether
\[
  \Hom_{\oO(n)}( \C[\R^{n\times k}],\pi)\ne 0?
\]
Note that
\[
  \C[\R^{n\times k}]=\bigotimes^k \C[\R^n].
\]
Since $\C[\R^n]$ contains the trivial representation of $\oO(n)$, we have
\be\label{ku0}
  \Hom_{\oO(n)}( \C[\R^{n\times k}], \pi)\neq 0\quad \Rightarrow\quad   \Hom_{\oO(n)}( \C[\R^{n\times (k+1)}], \pi)\neq 0.
\ee
Also we have a surjective homomorphism
\[
    \C[\R^{n\times n}]\xrightarrow{\textrm{restriction}}  \C[\oO(n)],
\]
and since $\pi$ always occurs in $ \C[\oO(n)]$, we see that
\be\label{ho0}
  \Hom_{\oO(n)}( \C[\R^{n\times k}], \pi)\neq 0\quad \textrm{if }k\geq n.
\ee
In view of \eqref{ku0} and \eqref{ho0}, we define the first occurrence index
\[
  \operatorname{n}(\pi):=\min\{k\geq 0\mid \Hom_{\oO(n)}( \C[\R^{n\times k}], \pi)\neq 0\}.
\]

The following result may be read off from the explicit description of the correspondence $\pi \leftrightarrow \sigma$ of Theorem \ref{sh2} (see \cite{KV}), and it is a special case of the conservation relations which we will discuss in the next section.

\begin{thm}
For any irreducible finite dimensional representation $\pi$ of $\oO(n)$, we have
\[
 \operatorname{n}(\pi)+\operatorname{n}(\pi\otimes \det)=n.
\]
\end{thm}

\begin{remark} The fact that $\operatorname{n}(\det)=n$ is implicit in Weyl's book \cite[Chapter 2]{Wey}.
\end{remark}

\section{Transcending classical invariant theory}

In this section, we give an informal introduction to Howe's paper of the same title \cite{Ho3}. We also describe two fundamental principles of the theory: the Howe duality theorem and Kudla-Rallis conservation relations.

We start from the setting of Section 1, where we have the representation
\begin{equation*}
 \oO(n)\times \s\p_{2k}(\C)\curvearrowright \C[\R^{n\times k}].
\end{equation*}
The above action of the symplectic Lie algebra $\s\p_{2k}(\C)$ does not integrate to a Lie group representation in general. This is to say that the space $\C[\R^{n\times k}]$ has local symmetry of the symplectic group, but no global symmetry. However, if we consider the space of Schwartz functions instead of the space of polynomials, we do have the global symmetry. We shall describe this in what follows.

In general, we may consider the real orthogonal group (definite or indefinite)
 \begin{eqnarray*}
     \oO(p,q)&:=&\left\{g\in \GL_{p+q}(\R)\mid g^t \begin{bmatrix}
                                                                                                 1_p & 0 \\
                                                                                                 0 & -1_q \\
                                                                                               \end{bmatrix} g= \begin{bmatrix}
                                                                                                1_p & 0 \\
                                                                                                 0 & -1_q \\
                                                                                               \end{bmatrix}\right\}, \quad (p,q\geq 0).
\end{eqnarray*}

Similar to \eqref{acto}, we have an obvious representation
 \be\label{actoon}
   \oO(p,q)\curvearrowright \CS(\R^{(p+q)\times k}),
\ee
where ``$\CS$" indicates the space of Schwartz functions.
In the same spirit as \eqref{actsp} but for more subtle reasons, the representation in \eqref{actoon} extends to a representation
 \begin{equation}\label{actspgg}
 \oO(n)\times \widetilde \Sp_{2k}(\R)\curvearrowright \CS(\R^{(p+q)\times k}),
\end{equation}
where $\widetilde \Sp_{2k}(\R)$ denotes the metaplectic double cover of the symplectic group $\Sp_{2k}(\R)$.

Let $\Irr(\oO(p,q))$ denote the set of isomorphism classes of irreducible representations (with certain technical conditions) of $\oO(p,q)$. Similar notations, such as $\Irr(\widetilde \Sp_{2k}(\R))$, will be used. The following result is the Howe duality theorem for the so-called dual pair $(\oO(p,q),\Sp_{2k}(\R))$, and we will discuss its general statement in the later part of this section.

\begin{thm}\label{hoos}
Denote
\begin{eqnarray*}
  \omega_k:&=&\mathcal S(\R^{(p+q)\times k})\\
  \Omega:&=&\{\pi\in  \Irr(\oO(p,q))\mid  \Hom_{\oO(p,q)}(\omega_k, \pi)\neq 0\},\\
  \Omega':&=&\{\pi'\in \Irr(\widetilde{\Sp}_{2k}(\R))\mid  \Hom_{ \widetilde{\Sp}_{2k}(\R)}(\omega_k, \pi')\neq 0\}.\\
\end{eqnarray*}
Then the relation
\[
  \{(\pi, \pi')\in \Irr(\oO(p,q))\times \Irr(\widetilde{\Sp}_{2k}(\R))\mid \Hom_{\oO(p,q)\times \widetilde{\Sp}_{2k}(\R)}(\omega_k, \pi\widehat \otimes \pi')\neq 0\}
\]
defines a one-to-one correspondence
\[
   \Irr(\oO(p,q))\supset \Omega\leftrightarrow \Omega'\subset \Irr(\widetilde{\Sp}_{2k}(\R)).
\]
Moreover, for any $(\pi, \pi')\in \Irr(\oO(p,q))\times \Irr(\widetilde{\Sp}_{2k}(\R))$, we have
\[
   \dim \Hom_{\oO(p,q)\times \widetilde{\Sp}_{2k}(\R)}(\omega_k, \pi\widehat \otimes \pi')\leq 1.
\]

\end{thm}

We proceed to describe another fundamental principle governing the correspondence of Theorem \ref{hoos}. Fix $\pi \in \Irr(\oO(p,q))$.
Similar to \eqref{ku0}, we have that
\be\label{ku1}
  \Hom_{\oO(p,q)}(\CS(\R^{(p+q)\times k}), \pi)\neq 0\quad \Rightarrow\quad   \Hom_{\oO(p,q)}( \CS(\R^{(p+q)\times (k+1)}), \pi)\neq 0.
\ee
This is called Kudla's persistence principle \cite{Ku1}. Similar to \eqref{ho0}, we have that
\be\label{ho1}
  \Hom_{\oO(p,q)}( \CS(\R^{(p+q)\times k}), \pi)\neq 0\quad \textrm{if }k\geq p+q.
\ee
The condition ``$k\geq p+q$" is called Howe's stable range condition.
In view of \eqref{ku1} and \eqref{ho1}, we define the first occurrence index
\begin{equation}\label{npi}
  \operatorname{n}(\pi):=\min\{k\geq 0\mid \Hom_{\oO(p,q)}(\CS(\R^{(p+q)\times k}), \pi)\neq 0\}.
\end{equation}

The following result is called conservation relations (for the group $\oO(p,q)$). It was conjectured by Kudla and Rallis \cite{KR}, and established by the authors in \cite{SZ}.

\begin{thm}
For any $\pi \in \Irr(\oO(p,q))$, we have
\[
 \operatorname{n}(\pi)+\operatorname{n}(\pi\otimes \det)=p+q.
\]
\end{thm}

\begin{remark} The paper \cite{SZ} proves a conservation relation for a (type I) classical group $G$, which is valid for any irreducible smooth representation of $G$, and for any local field of characteristic zero.
\end{remark}

Now we come to the general statement of Howe duality conjecture \cite{Ho1}. Let $W$ be a finite-dimensional symplectic vector space over $\R$ with symplectic form $\langle \,,\, \rangle_W$.
Denote by $\tau$ the involution of $\End_\R(W)$ specified by
\[
   \langle x\cdot u, v\rangle_W=\langle u, x^\tau\cdot v\rangle_W, \qquad u,v\in W,\,x\in \End_\rF(W).
\]
Let $(A, A')$ be a pair of $\tau$-stable semisimple $\R$-subalgebras of $\End_\R(W)$ which are mutual centralizers of each other. Put $G:=A\cap \Sp(W)$ and $G':=A'\cap \Sp(W)$, which are closed subgroups of $\Sp(W)$. Following Howe, we call the pair $(G,G')$ so obtained a {\em reductive dual pair} in $\Sp(W)$.
We say that the pair $(A, A')$ (or the reductive dual pair $(G,G')$) is irreducible of type I if  $A$ (or equivalently $A'$) is a simple algebra, and that it is  irreducible of type II if  $A$ (or equivalently $A'$) is the product of two simple algebras which are exchanged by $\tau$. A complete classification of such dual pairs was also given by Howe \cite{Ho1}.

Denote by
$\mathrm H(W):=W\times \R$ the Heisenberg group attached to $W$, with
group multiplication
\[
      (u, \alpha)(v,\beta):=(u+v, \alpha +\beta+\langle u,v\rangle_W), \qquad u,v\in W, \ \alpha,\beta\in \R.
\]
Its center is identified with $\R$ in the obvious way.

Fix an arbitrary non-trivial unitary character $\psi: \R\rightarrow \mathbb C^\times$.
\begin{thm}\label{vn}(Stone-von Neumann)
Up to isomorphism, there is a unique irreducible
unitary representation $\widehat \omega_\psi$ of $\mathrm H(W)$ with the central character $\psi$.

\end{thm}

Let $\widehat \omega_\psi$ be as in Theorem \ref{vn}. Set
\[
   \omega_\psi:=\{v\in \widehat \omega_\psi \mid \textrm{the map $\rH (W)\rightarrow \widehat \omega_\psi,\ g\mapsto g\cdot v$ is smooth}\}.
\]
This is a smooth representation of $\rH(W)$ with the usual smooth topology.

Let $(G,G')$ be a reductive dual pair in $\Sp(W)$. Let $\widetilde G$ and $\widetilde G'$ be finite fold covering groups of $G$ and $G'$, respectively. Then $\widetilde G\times \widetilde G'$ acts on $\rH(W)$ as group automorphisms via the action of $G\times G'$ on $W$. Using this action, we form the semidirect product
\[
  (\widetilde G\times \widetilde G')\ltimes \rH(W).
\]
We say that a representation of $(\widetilde G\times \widetilde G')\ltimes \rH(W)$ is a smooth oscillator representation if its restriction to $\rH(W)$ is isomorphic to  $\omega_\psi$.
We say that a representation of $\widetilde G\times \widetilde G'$ is a smooth oscillator representation if it extends to a smooth oscillator representation of $(\widetilde G\times \widetilde G')\ltimes \rH(W)$. Basic facts about smooth oscillator representations may be found in \cite{Wei,Ku2}.

The following is the fundamental result of Howe \cite{Ho2}, establishing his famous conjecture for all real reductive dual pairs.

\begin{thm}\label{hod}
Let $\omega$ be a  smooth oscillator representation of  $(\widetilde G\times \widetilde G')\ltimes \rH(W)$.
Denote
\begin{eqnarray*}
  \Omega:&=&\{\pi\in  \Irr(\widetilde G)\mid  \Hom_{\widetilde G}(\omega, \pi)\neq 0\},\\
  \Omega':&=&\{\pi'\in  \Irr(\widetilde G')\mid  \Hom_{\widetilde G'}(\omega, \pi')\neq 0\}.\\
\end{eqnarray*}
Then the relation
\[
  \{(\pi, \pi')\in  \Irr(\widetilde G)\times  \Irr(\widetilde G'))\mid \Hom_{\widetilde G\times \widetilde G'}(\omega, \pi\widehat \otimes \pi')\neq 0\}
\]
defines a one-to-one correspondence
\[
   \Irr(\widetilde G)\supset \Omega\leftrightarrow \Omega'\subset \Irr(\widetilde G').
\]
Moreover, for any $(\pi, \pi')\in \Irr(\widetilde G)\times \Irr(\widetilde G')$, we have
\[
   \dim \Hom_{\widetilde G\times \widetilde G'}(\omega, \pi\widehat \otimes \pi')\leq 1.
\]
\end{thm}

The correspondence specified in the above theorem is called the local theta correspondence, or theta correspondence in short.

\begin{example}
The groups $\oO(p,q)$ and $\Sp_{2k}(\R)$ form a reductive dual pair in $\Sp_{2k(p+q)}(\R)$. When $\widetilde G=\oO(p,q)$ and $\widetilde G'=\widetilde \Sp_{2k}(\R)$, Theorem \ref{hod} reduces to Theorem \ref{hoos}.

\end{example}

\begin{example}
The groups $\GL_m(\R)$ and $\GL_{n}(\R)$ form a reductive dual pair in $\Sp_{2mn}(\R)$. The following representation is a smooth oscillator  representation
 \[
   \GL_m(\R)\times \GL_n(\R)\curvearrowright \CS(\R^{m\times n}), \quad ((g,h)\cdot f)(x):=f(g^{-1} x h),
 \]
 where $g\in \GL_m(\R)$, $h\in \GL_n(\R)$, $x\in \R^{m\times n}$, $f\in  \CS(\R^{m\times n})$. Suppose $m\leq n$. Then the local theta correspondence yields an injective map
 \[
   \Irr(\GL_m(\R))\hookrightarrow \Irr(\GL_n(\R)).
 \]
 If $m=n$, then this is the bijection which sends $\pi\in \Irr(\GL_m(\R))$ to its contragradient representation $\pi^\vee$.

\end{example}

\begin{remark} Howe made his duality conjecture for both real and $p$-adic local fields in \cite{Ho1}. For the latter, the Howe duality conjecture has also been established completely, thanks to the works of Waldspurger \cite{Wa}, Minguez \cite{Mi}, Gan and Takeda \cite{GT}, and Gan and Sun \cite{GS}. The companion statement of multiplicity one was due to Waldspurger \cite{Wa} and Li, Sun and Tian \cite{LST}. We refer the readers to \cite{MVW} on basic approaches to local theta correspondence for $p$-adic local fields.
\end{remark}

\section{The explicit correspondence}

For applications to automorphic forms and unitary representation theory, it is of great interest to explicitly describe the local theta correspondence, in terms of Langlands-Vogan parametrization of irreducible representations \cite{Vo}.

The first important cases concern the so-called (irreducible) compact dual pairs, namely
\[(G,G')=(\oO(p), \Sp(2n, \R)), (\oU(p), \oU(r,s)), (\Sp(p), \oO^*(2n)).\]
The representations of $\cover{G'}$ arising from the correspondence are the so-called unitary lowest weight modules. The explicit correspondence is described in two landmark papers on the subject of unitary lowest weight modules, one by Kashiwara and Vergne \cite{KV} and another by Enright, Howe and Wallach \cite{EHW}.

For non-compact dual pairs, the situation is considerably more complicated. (The case of compact dual pairs is almost always a stepping stone for the general case.) In the following, we list some of the most relevant works, which give complete description of the correspondence for specific dual pairs (or specific classes of representations).
\begin{itemize}
\item Dual pairs $(G,G')$ of type I, where the ``size" of $G$ is not greater than that of $G'$: for a ``sufficiently regular" (genuine) discrete series representation $\pi$ of $\cover{G}$, Li shows that $\pi$ occurs in the correspondence, and the corresponding representation $\pi '$ of $\cover{G'}$ is a (explicitly described) unitary representation with nonzero cohomology. See \cite{Li} for details.
\item Dual pairs of type II: complete correspondence is given by Moeglin \cite{Mo} (over the real), Adams-Barbasch \cite{AB1} (over the complex) and Li-Paul-Tan-Zhu \cite{LPTZ} (over the quaternion).
\item Complex dual pairs of type I: complete correspondence is given by Adams-Barbasch \cite{AB1}.
\item The dual pair $(\oO(p,q), \Sp(2n,\R))$ in the so-called ``(almost) equal rank" case: complete correspondence is given by Moeglin and Paul \cite{Mo,Pa2} (for $p+q=2n, 2n+2$), and  Adams-Barbasch \cite{AB2} (for $p+q=2n+1$).
\item The dual pair $(\oU(p,q),\oU(r,s))$ in the ``equal rank" case: complete correspondence is given by Paul \cite{Pa1}, where $p+q=r+s$.
\item The dual pair $(\Sp(p,q),\oO^*(2n))$ in the ``equal rank" case and beyond: complete correspondence is given by Li-Paul-Tan-Zhu \cite{LPTZ}, where $p+q\leq n$.
\item Other cases: complete correspondence is given by Przebinda \cite{Pr} for the dual pair $(\oO(2,2), \Sp(4,\R))$, Bao \cite{Ba} for the dual pair $(\Sp(p,q),\oO^*(4))$  and Fan \cite{Fa} for the dual pair $(\oO(p,q), \Sp(2n,\R))$ (where $p+q=4$).
\end{itemize}

\begin{remark} Common techniques in determining the explicit correspondence include (a) correspondence of discrete series as in the work of Li \cite{Li}; (b) the induction principle of Kudla \cite{Ku1}; and (c) the correspondence of $\cover{K}$ and $\cover{K'}$-types in the space of joint harmonics \cite{Ho3}, and arguments to identify representations using minimal $K$-types.
\end{remark}

\section{Automatic continuity}
In the literature, the theta correspondence is mostly studied in the algebraic setting of Harish-Chandra modules. In fact Howe proved his famous duality theorem in the algebraic setting which he shows to imply the duality theorem in the smooth setting \cite{Ho2}. Ever since this, an ``automatic continuity" result has been expected, namely that the algebraic version of theta correspondence agrees with the smooth version of theta correspondence.

We shall highlight this ``automatic continuity" problem. Our exposition follows that of \cite{BS}. We start with the following

\begin{prp}\label{carin}
Let $(G, G')$ be a reductive dual pair in the real symplectic group $\Sp(W)$.  Up to conjugation by  $G\times G'$, there exists a unique Cartan involution $\theta$ of $\mathrm{Sp}(W)$  such that
\begin{equation}\label{Cartan Decomposition 2}
    \theta(G)=G\quad\textrm{and}\quad\theta(G')=G'.
\end{equation}
\end{prp}

Now let $\theta$ be as in Proposition \ref{carin}.  Then $\theta$ restricts to Cartan involutions on $G$ and $G'$.  Denote the fixed point groups by
\[
  C:=(\mathrm{Sp}(W))^\theta,  \quad K:=G^\theta \quad \textrm{and}\quad  K':=(G')^\theta .
\]
They are maximal compact subgroups of $\mathrm{Sp}(W)$, $G$ and $G'$, respectively. Write
\[
\textrm{$\widetilde K\rightarrow K$, $\quad \widetilde K'\rightarrow K'\quad $ and $\quad \widetilde C\rightarrow C$}
 \]
 for the covering maps induced by the covering maps
 \[
 \textrm{$\widetilde G\rightarrow G$, $\quad \widetilde G'\rightarrow G'\quad $ and $\quad\widetilde \Sp(W)\rightarrow \Sp(W)$,}
 \]
  respectively.

As in Theorem \ref{hod}, let $\omega$ be a  smooth oscillator representation of  $(\widetilde G\times \widetilde G')\ltimes \rH(W)$.
Denote
\begin{eqnarray*}
  \Omega&:=&\{\pi\in  \Irr(\widetilde G)\mid  \Hom_{\widetilde G}(\omega, \pi)\neq 0\},\\
  \Omega'&:=&\{\pi'\in  \Irr(\widetilde G')\mid  \Hom_{\widetilde G'}(\omega, \pi')\neq 0\},\\
 \CR&:=& \{(\pi, \pi')\in  \Irr(\widetilde G)\times  \Irr(\widetilde G')\mid \Hom_{\widetilde G\times \widetilde G'}(\omega, \pi\widehat \otimes \pi')\neq 0\}.
\end{eqnarray*}
It is known \cite{Wei} that the representation $\omega|_{\rH(W)}$ uniquely extends to a representation of $\widetilde \Sp(W)\ltimes \rH(W)$.  Write $\omega^{\mathrm{alg}}$ for the space of $\widetilde C$-finite vectors in $\omega$.
For any $\pi\in \mathrm{Irr}(\widetilde{G})$, define
\begin{equation}\label{homalg}
\mathrm{Hom}_{\widetilde{G}}^{\mathrm{alg}}(\omega,\pi):=\mathrm{Hom}_{\mathfrak g,\widetilde{K}}(\omega^{\mathrm{alg}},\pi^\mathrm{alg}),
\end{equation}
where $\mathfrak g$ denotes the complexified  Lie algebra of $G$,  and $\pi^\mathrm{alg}$ denotes the $(\mathfrak g, \widetilde{K})$-module of $\widetilde{K}$-finite vectors in $\pi$. We define $\mathrm{Hom}_{\widetilde{G'}}^{\mathrm{alg}}(\omega,\pi')$ and $\mathrm{Hom}_{\widetilde{G}\times \widetilde{G'}}^{\mathrm{alg}}(\omega,\pi\widehat{\otimes} \pi')$ in analogous ways, where $\pi'\in \mathrm{Irr}(\widetilde{G'})$.
In this algebraic setting, we let
\begin{eqnarray*}
  \Omega^{\mathrm{alg}} &:=& \{\pi\in \mathrm{Irr}(\widetilde{G})\mid \mathrm{Hom}^{\mathrm{alg}}_{\widetilde{G}}(\omega,\pi)\neq 0\}, \\
  \Omega'^{\mathrm{alg}}&:=& \{\pi'\in \mathrm{Irr}(\widetilde{G'})\mid \mathrm{Hom}^{\mathrm{alg}}_{\widetilde{G'}}(\omega,\pi')\neq 0\}, \\
 \CR^{\mathrm{alg}}&:= & \{(\pi,\pi')\in \mathrm{Irr}(\widetilde{G})\times \mathrm{Irr}(\widetilde{G'})\mid \mathrm{Hom}^{\mathrm{alg}}_{\widetilde{G}\times \widetilde{G'}}(\omega,\pi\widehat \otimes \pi')\neq 0\}.
\end{eqnarray*}
It is clear that
\begin{equation}\label{inclusion}
\left\{
  \begin{array}{l}
        \Omega\subset \Omega^{\mathrm{alg}},\\
    \Omega'\subset \Omega'^{\mathrm{alg}},\\
    \CR\subset \CR^{\mathrm{alg}}.\\
  \end{array}
\right.
\end{equation}

As a prelude of Theorem \ref{hod}, Howe proved the following algebraic analog. As in \cite[Section 2]{Ho3}, this implies Theorem \ref{hod}.

\begin{thm}\label{hod2}
 The relation $\CR^{\mathrm{alg}}$
defines a one-to-one correspondence
\[
   \Irr(\widetilde G)\supset \Omega^{\mathrm{alg}}\leftrightarrow \Omega'^{\mathrm{alg}}\subset \Irr(\widetilde G').
\]
Moreover, for any $(\pi, \pi')\in \Irr(\widetilde G)\times \Irr(\widetilde G')$, we have
\[
   \dim \Hom^{\mathrm{alg}}_{\widetilde G\times \widetilde G'}(\omega, \pi\widehat \otimes \pi')\leq 1.
\]

\end{thm}

We state the following folklore conjecture on the coincidence of the smooth version and the algebraic version of theta correspondence, which has been expected since the publication of \cite{Ho3}.

\begin{conj}\label{Conjecture on Coincidence0}
The inclusions in \eqref{inclusion} are all equalities.
\end{conj}

Recall from Section 2 the construction of a reductive dual pair $(G,G')\subset \Sp(W)$ from a mutually centralizing pair $(A, A')$ of $\tau$-stable semisimple $\R$-subalgebras of $\End_\R(W)$, as well as the notion of type I and type II dual pairs. Let $\mathbb H$ denote a quaternionic division algebra over $\mathbb R$, which is unique up to isomorphism.

The following theorem is the main result of \cite{BS}.

\begin{thm}\label{Main Theorem}
Conjecture $\ref{Conjecture on Coincidence0}$ holds when the dual pair $(G, G')$ contains no quaternionic type I irreducible factor, that is, when $A$ contains no $\tau$-stable ideal which is isomorphic to a matrix algebra $\mathrm M_{n\times n}(\mathbb H)$ $($$n\geq 1$$)$.
\end{thm}

\begin{remark} For the quaternionic dual pair $(G, G')=(\Sp (p,q), \oO^*(2n))$, it is also shown in \cite{BS} that Conjecture $\ref{Conjecture on Coincidence0}$ holds if $p+q\leq n$.
\end{remark}

We note that the main tool of \cite{BS} showing the conjecture on coincidence is the conservation relations established in \cite{SZ}. In the following we sketch a proof for orthogonal-symplectic dual pairs.

Let $G=\oO(p,q)$ and $\pi\in \Irr(G)$. Define the first occurrence index $\mathrm n(\pi)$ as in \eqref{npi}, and its algebraic analogue
 \[
  \mathrm{n}'(\pi):= \min\left\{n\geq 0\mid   \mathrm{Hom}_{G}^\mathrm{alg}(\mathcal S(\mathbb R^{(p+q)\times n}),\pi)\neq 0\right\},
\]
where $\mathrm{Hom}_{G}^\mathrm{alg}$ is defined as in \eqref{homalg}, using an obvious subspace $\mathcal S(\mathbb R^{(p+q)\times n})^\mathrm{alg}$ of $\mathcal S(\mathbb R^{(p+q)\times n})$.

It is clear that $ \mathrm{n}'(\pi)\leq \mathrm n(\pi)$. Similar to the conservation relation
\[
 \operatorname{n}(\pi)+\operatorname{n}(\pi\otimes \det)=p+q,
\]
we may prove the algebraic analog:
\[
 \operatorname{n}'(\pi)+\operatorname{n}'(\pi\otimes \det)=p+q.
\]
Thus $\mathrm n(\pi)=\mathrm{n}'(\pi)$. This easily implies that Conjecture $\ref{Conjecture on Coincidence0}$ holds for orthogonal-symplectic dual pairs.

\end{document}